\documentclass{amsart}
\usepackage{latexsym,amssymb}
\usepackage{graphicx}
\usepackage{epsfig}
\input epsf
\usepackage{amsfonts}
\usepackage{amscd}
\usepackage{amsmath}

\theoremstyle{plain} \numberwithin{equation}{section}
\newtheorem{thm}{Theorem}[section]

\newtheorem{lem}[thm]{Lemma}

\theoremstyle{definition}
\newtheorem{de}[thm]{Definition}

\newtheorem{rem}[thm]{Remark}
\newtheorem{Example}[thm]{Example}

\newcommand{\R}{\mathbb{R}}

\newcommand{\Z}{\mathbb{Z}}

\newcommand{\cA}{\mathcal{A}}

\newcommand{\cD}{\mathcal{D}}
\newcommand{\cF}{\mathcal{F}}
\newcommand{\cG}{\mathcal{G}}
\newcommand{\cB}{\mathcal{B}}

\newcommand{\cM}{\mathcal{M}}
\newcommand{\cN}{\mathcal{N}}

\newcommand{\cT}{\mathcal{T}}

\newenvironment{deflist}
{\begin{list}{$\blacktriangleright$}{%
\setlength{\leftmargin}{17pt} 
\setlength{\labelsep}{6.5pt}
\setlength{\labelwidth}{10pt}
\setlength{\parsep}{2pt}
\setlength{\topsep}{0.2\baselineskip}
\setlength{\itemsep}{0.2\baselineskip}}}
{\end{list}}

\newcommand{\GL}{\mathrm{GL}}

\newcommand{\id}{\mathrm{id}}
\newcommand{\supp}{\mathrm{supp}}

\begin{document}

\begin{abstract}
In this short note we discuss the interplay between finite Coxeter groups
and construction of wavelet sets, generalized multiresolution analysis
and sampling.
\end{abstract}

\author{Mihaela Dobrescu}
\email{mihaela.dobrescu@furman.edu}
\address{Department of Mathematics\\
Christopher Newport University\\
Newport News, VA 23606, USA
}
\author{Gestur \'Olafsson}
\email{olafsson@math.lsu.edu}
\address{Department of Mathematics\\
Louisiana State University\\
Baton Rouge, LA 70803, USA}

\title{Coxeter Groups, Wavelets, Multiresolution and Sampling}

\thanks{M. Dobrescu was partially
supported by DMS-0139783.
The research of G. \'Olafsson was supported by NSF grants DMS-0139783 and DMS-0402068.}
\maketitle
\section*{Introduction}
\noindent
Finite reflection groups are examples of finite Coxeter groups.
Those groups show up in a natural way in geometry as symmetry groups of
geometric objects, fractal geometry in the classification of simple Lie algebras,
in representation theory, theory of special functions and
other places in analysis, cf. \cite{PG06,Mass94,OP04,op06,Op88,Op93}
and the reference therein for few examples. On the other
hand, the only example we know of, where those groups have
shown up in the connection with wavelet theory and
multiresolution analysis is in the book by P. Massopust \cite{Mass94} and
related construction by him and his coworkers.

Let $\mathcal {D}\subseteq GL(n,\mathbb{R})$ and
$\mathcal{T}\subseteq\mathbb{R}^n $ countable sets.
Recall, that a $(\mathcal D,\mathcal T)$-\textit{wavelet} is a square integrable
function $\psi$ with  the property that the set
\begin{equation}\label{eq:wavelet}
\{|\det d|^{\frac{1}{2}}\varphi(dx+t)\mid d\in\mathcal{D},t\in\mathcal{T}
\}
\end{equation}
forms an orthonormal basis for $L^{2}(\mathbb{R}^{n})$.
A special class of wavelets are the ones corresponding to
\textit{wavelet sets}. Those are functions $\psi$ such that
$\cF (\psi)=\chi_{\Omega}$ for a measurable subset $\Omega$ of $\R^n$.
The set $\mathcal {D}$ is then called \textit{the dilation set} and the set
$\mathcal{T}$ is called \textit{the translation set}. Quite often, one
assume that the
dilation set is a group and that the translation set is a latices,
i.e., a discrete subgroup $\Gamma\subset \R^n$ such
that $\R^n/\Gamma$ is compact. In particular, the
simplest example is the group generated by
one element, $\cD=\{a^k\mid k\in \Z\}$.
In \cite{O02,olafsson_speegle:03} more general sets of dilations were
considered, and in general those dilations do not form a group.
Even more general constructions can be found in \cite{A03}.

In this article we consider the case where $\cD$ is of
the form $\cD=\{a^k\id \mid k\in\Z\}W$, where
$a$ is an expansive matrix and $W$ is a finite Coxeter group,
see Section \ref{section_1} for the definition. We use results from
our previous article \cite{DO2005}. Those results are discussed in
Section \ref{section_2}. We would like to remark, that this construction
is more general than needed here. In particular it is not
needed that $\cD$ is a group.

Most of the examples of wavelet sets tend to be fractal like and symmetric around $0$.
Our aim is to construct wavelet functions that some directional properties in
the frequency domain. The finite Coxeter group is then
used to rotate the frequency domain to cover all of $\R^n$.
Some two-dimensional examples are discussed at
the end of Section \ref{section_2} to motivate the construction
in later sections.
The construction is still fractal, but it is not symmetric around $0$ anymore.
The generalization of a rotation group for dimensions higher than two is a Coxeter group.
The construction is generalized to higher dimensions in Section \ref{section_1}
and  Section \ref{section_5}.

A natural question one asks  when working with wavelets is if they are related
with any multiresolution analysis.
It is well known that one needs multiwavelets when working in higher dimensions which
means that the wavelets coming from
the wavelets sets described above are not associated with any multiresolution analysis.

In Section \ref{section_3} and Section \ref{section_4} we discuss the
construction of scaling sets and the
associated multiresolution analysis and multiwavelets for
this situation. In particular, our wavelet is still directional in
the frequency domain. In fact the support
of the Fourier transform is supported in cones which are
fundamental domains for the action of a Coxeter group on the Euclidian space.
We are also illustrating our results with a few examples.

In the final section we discuss the relation between these results and sampling theory.
Any square integrable function can be written as a sum of its projections on subspaces,
where each subspace contains only signals supported in the frequency
domain in the cones mentioned above.
Each projection can then be sampled using a version of the
Whittaker-Shannon-Kotel'nikov sampling theorem \cite{Sh48,Za93}
for spectral sets stated in Theorem \ref{wsk}.

\section{Existence of subspace wavelet sets}\label{section_2}
\noindent
In this section we recall some
general results from \cite{DO2005} which are the basic for
the construction later in this work. Most of the literature
deals with dilations groups generated by one element. The main
idea here is to consider dilation sets that can be factorized
as a product of finitely many groups (or more general set)
and then use inductive construction
to reduce the general case to the simple one.
Let us remark that the statements in this section hold for
more general settings, i.e., one could replace $\R^n$ by a measure space $M$
and $\GL (n,\R)$ by a group of automorphisms of $M$.

For $\cA,\cB\subset \GL (n,\R)$ we say that the product
$\cA\cB=\{ab\mid a\in \cA,\, b\in\cB\}$ is \textit{direct} if
$a_1b_1=a_2b_2$, $a_1,a_2\in\cA$, $b_1,b_2\in\cB$, implies that
$a_1=a_2$ and $b_1=b_2$. For the proof of the following
statements see \cite{DO2005}.

\begin{de}
A measurable tiling of a measure space $(M,\mu)$ is a countable collection of
subsets $\{\Omega_j\}$ of $M$, such that
\begin{displaymath}
\mu(\Omega_i\cap\Omega_j)=0\, ,
\end{displaymath}
for $i\not= j$,
and
$$
\mu(M\diagdown\bigcup_{j}\Omega_j)=0\, .$$
If $\Omega\subset \R^n$ is measurable and $\cA$ a set of
diffeomorphism of $\R^n$, then
$\Omega$ is a $\cA$-tile if $\{d (\Omega) \mid d\in \cA\}$ is a measurable
tiling of $\R^n$.
\end {de}

\begin{lem}\label{th-2.1}
Let $M\subseteq \R^n$ be measurable.
Let $\mathcal{A},\mathcal{B}\subset \GL(n,\mathbb{R})$ be two non-empty sets, such that
the product $\cA\cB$ is direct.
Let $\mathcal{D} = \mathcal{A}\mathcal{B}=\{ab\mid a\in
\mathcal{A}\, ,\,\, b\in \mathcal{B}\}$.
Then there exists a $\mathcal{D}$-tile $\Omega$ for $M$ if and only if there
exists a measurable set $N\subseteq \R^n$,
such that $\mathcal{A} N$ is a measurable tiling of
$M$, and a
$\mathcal{B}$-tile $\Omega$ for $N$.
\end{lem}

\begin{rem} We would like to remark at this point, that we do
not assume that $\Omega\subseteq M$, nor that $N\subseteq M$.
But this will in fact be the case in most applications because
$\cD$ will contain the identity matrix.
\end{rem}

\begin{thm}[(Construction of wavelet sets by steps, I)]\label{th-2.2}
Let $\mathcal{M},\mathcal{N}\subset \GL(n,\R)$ be two non-empty
subsets such that the product $\cM\cN$ is direct.
Let $\mathcal{L} = \mathcal{M}\mathcal{N}$.
Assume that $M\subseteq \mathbb{R}^n$ with $|M|>0$, is measurable.
Let $\mathcal{T}\subset \mathbb{R}^n$ be discrete.
Then there exists a $(\mathcal{L} ,\mathcal{T})$-wavelet
set $\Omega\subset M$ for $M$ if and only if there
exists a $\mathcal{N}^T$-tiling set $N\subset M$ and a
$(\mathcal{M},\mathcal{T})$-wavelet set $\Omega_1$ for $N$.
\end{thm}

Recall that if $\cD\subseteq \GL(n,\R)$ and $\cG\subset \GL (n,\R)$ is
a group that acts on $\cD$ from the right, then there exists
a subset $\cD_1\subseteq \cD$, such that $\cD=\cD_1\cG$ and the
product is direct. Note that we do not assume that $\cG\subset \cD$.

\begin{thm}[(Construction of wavelet sets by steps, II)] Let $\cD\subset
\GL (n,\R)$ and $M\subseteq \R^n$ measurable
with $|M|>0$. Let $\cT\subset \R^n$ be discrete.
Assume that $\mathcal{G}\subset \GL(n,\R)$ is a group that
acts on $\cD$ form the right. Let $\cD_1\subseteq \cD$ be
such that $\cD=\cD_1\cG$ as a direct product.
Then there exists a $(\cD,\cT)$-wavelet set
$\Omega$ for $M$ if only only if there exists
a $\cG^T$-tiling set $N$ for $M$ and a $(\cD_1,\cT)$-wavelet set $\Omega_1$ for
$N$.
\end{thm}

For a measurable set $M\subseteq \R^n$ such that
$\overline{M}\setminus M$ has measure zero, let
$$L_M^2(\R^n)=\{f\in L^2(\R^n)\mid \supp (\cF(f))\subseteq M\}\, .$$
The question is then, how to obtain a wavelet set for
the starting subset $N$. The following result gives one
way to do that.

\begin{thm}[(Existence of subspace wavelet sets)]\label{th-exwave}
Let $M \subseteq \mathbb{R}^{n}$ be a measurable
set, $|M|>0$. Let $a\in \GL(n,\R)$ be an expansive
matrix and $\emptyset\not=\mathcal{D}\subset \GL (n,\mathbb{R}) $.
Assume that
$\cD^T$ is  a multiplicative tiling of $M$, $a\cD=\cD$ and
$a^T M=M$. If $\mathcal{T}$ is a lattice,
then there exists a
measurable set $\Omega \subseteq M$ such that $ \Omega + \cT$ is a
measurable tiling
of $\mathbb{R}^{n}$
and $\cD^T \Omega$  is a measurable tiling of
$M$. In particular, $\Omega$ is a $L^2_M(\R^n)$-subspace $(\cD,\cT)$-wavelet set.
\end{thm}

\begin{Example}\label{Ex-3.1}
We would like to note here that in general, for a given set $\cD$, there are
several ways to decompose it as a direct product. As an example take
the set - which in fact is a group -
$$\cD=\cD_{a,m} = \{a^kR_{2 \pi j/m }\mid k\in\Z, j=0,\ldots m-1\}\, .$$
Here $a>1$ and $R_{\theta }$ stands for the rotation
$$R_{\theta }=\begin{pmatrix} \cos \theta & -\sin \theta\\
\sin \theta  &\cos \theta \end{pmatrix}\, .$$
We can take $\cG=\{R_{2 \pi j/m }\mid j=0,\ldots , m-1\}$
and $\cD_1=\{a^k\id \mid k\in \Z\}$. Then we can take
\begin{eqnarray*}
N=\R^2_{2\pi /m}&=&\{r(\cos \theta ,\sin \theta)^T\mid 0 \le  \theta \le 2\pi /m,\, r>0 \}\\
&=&\{(x,y)\in \R^2\mid 0\le  x,y\quad \text{and}\quad y\le x\tan (2\pi /m)\}\, .
\end{eqnarray*}
But we could also take $\cG$ as above and
$$\cD_1=\{(aR_{2\pi j/m})^k \mid k\in \Z\}$$
for some $0\le j\le m-1$. In this case we would take
$$N=\bigcup_{k\in \Z}(aR_{2\pi j/m})^k\{r(\cos \theta ,\sin \theta )^T\mid
a\le r\le a^2\quad \text{and}\quad 0\le \theta \le 2\pi /m\}$$
whose interior is not connected.
Note that if we take $a=\sqrt{2}$ and $m=4$, so $R_{2\pi /m}=
R_{\pi /4}$, then
$$aR_{\pi /4}=\begin{pmatrix} 1 & -1 \cr 1 & 1\end{pmatrix}$$
which often shows up in examples. As $a>1$ it follows that
$a\id $ and $aR_{2\pi j/m}$ are expansive matrices, Theorem \ref{th-exwave}
implies that in both cases and for any full rank lattice $\Gamma$, the
set $N$ contains a wavelet set $\Omega$, which will be quite different for the two cases.

Let $\mathcal{T}=\mathbb{Z}^2$, $\cD=\cD_{a,m} = \{a^kR_{2 \pi j/m }
\mid k\in\Z, j=0,\ldots m-1\}\, $, and let
\begin{displaymath}
\begin{array}{rcl}
E & = & [0,1]\times[0,\tan({2\pi /m})]\,\,\,\, if \,\,\,m\neq 2,4\\[1ex]
E & = & [0,1]^2 \,\,\,\, if \,\,\,m=4\\[1ex]
E & = &[-1,1]\times[0,1]\,\,\,\, if \,\,\,m=2\\[1ex]
F & = & \{(x,y)\in \R^2_{2\pi /m}|1<x<a\}.
\end{array}
\end{displaymath}
The wavelet set $\Omega$ has the form
\begin{displaymath}
\Omega=\bigcup_{i=1}^{2}\bigcup_{j=1}^{\infty}\Omega_{i,j} ,
\end{displaymath}
see figure  \ref{wavrot}.
The description of the $\Omega_{i,j}$ is as follows
\begin{displaymath}
\begin{array}{rcl}
\Omega_{1,1} & = & (E\setminus a^{-1}E)+(1,0)\\[1ex]
\Omega_{2,1} & = & a^{-2}\Big(F\setminus\big(E+(0,1)\big)\Big)\\[1ex]
\Omega_{1,2} &= & [(a^{-1}E\setminus a^{-2}E)\setminus \Omega_{2,1}]+(1,0)\\[1ex]
\Omega_{2,2} &= & a^{-3}[\Omega_{2,1}+(0,1)].
\end{array}
\end{displaymath}
For $j\geq 3$, we have the following formulas
\begin{displaymath}
\Omega_{1,j}=[(a^{-n+1}E\setminus a^{-n}E)\setminus \Omega_{2,j-1}]+(1,0)
\end{displaymath}
and
\begin{displaymath}
\Omega_{2,j}=a^{-n-1}[\Omega_{2,j-1}+(0,1)].
\end{displaymath}
From the construction, it is clear that $\Omega$ and $E$ are
$\mathcal{T}$-\textit{translation congruent},
and $\Omega$ and $F$ are $\mathcal{D}$-\textit{dilation congruent}.
 On the other hand, $F$ is a
$\mathcal{D}_{a,m}$-multiplicative tile and $\{E,\mathcal{T} \}$ is a spectral pair.
It follows that $\Omega$ is a $\mathcal{D}_{a,m}$-multiplicative tile
and $\{\Omega,\mathcal{T}\}$ is a spectral pair.
Thus $\Omega$ is a $(\mathcal{D}_{a,m}, \mathcal{T})$
wavelet set.

\begin{figure}[!ht]
\begin{center}
\epsfig{scale=0.7,angle=0,file=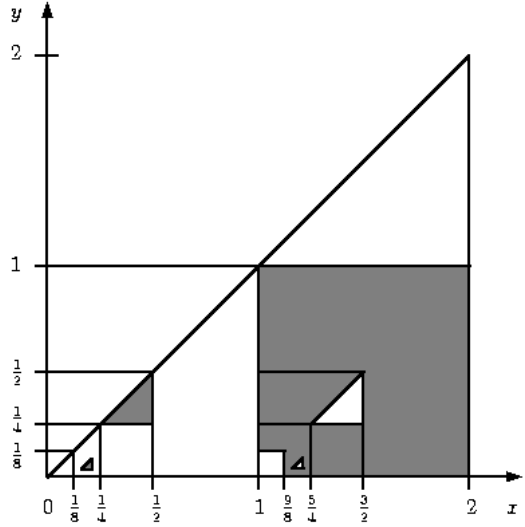}
\caption{A $(\mathcal{D}_{2,4},\mathcal{T})$ wavelet set.}
\label{wavrot}
\end{center}
\end{figure}

\end{Example}

\section{Coxeter Groups}\label{section_1}
\noindent
In the last example the group $\cG$ was a special case of a Coxeter group,
which we will introduce in  this section. The generalization to higher
dimensions of Example \ref{Ex-3.1} are the sets of the form
$\cD =\{a^kw\mid k\in \Z,\, w\in W(\Delta)\}$. The set $N$ is then
a convex cone which is a fundamental domain for $W$. Note that this
decomposition of $\cD$ is not the only one and we could also choose
other sets for $N$. We could also replace $a$ by any expansive matrix
$a$ such that $aN=N$.
We start
by collecting some well known facts on root systems and Coxeter groups.
We use \cite{Humphreys:92} as standard reference, but would also like to point out
\cite{bourLie4-6,He78,Mass94}.

A finitely generated group $W$
defined by the relations $(r_i r_j)^{m_{ij}}=1$, where $m_{ij}\in \mathbb{N}$,
$m_{ii}=1$ and $m_{ij}=m_{ji}$, is called a \textit{Coxeter group}. The finite
Coxeter groups can be realized as finite reflexion groups in $\mathrm{O}(n)$.
Let $\alpha \in \mathbb{R}^n$, $\alpha \not=0$. A \textit{reflection} along
$\alpha$ is a linear map  $r: \mathbb{R}^n\to \mathbb{R}^n$ such
that
\begin{itemize}
\item[$\bullet$] $r(\alpha )=-\alpha$;
\item[$\bullet$] The space $\{\lambda \in \mathbb{R}^n\mid r(\lambda )=\lambda \}$ is
a hyperplane in $\mathbb{R}^n$.
\end{itemize}
Note that a reflection $r$ is a non-trivial element of $\mathrm{O}(n)$ of order $2$.
The above reflection is given by
$$ \lambda \mapsto \lambda - \frac{2 (\lambda ,\alpha )}{(\alpha ,\alpha )}\alpha =:
r_\alpha (\lambda)\, .$$
Furthermore,
$$\{\lambda \in \mathbb{R}^n\mid r(\lambda )=\lambda \}=H_\alpha:=\{\lambda\in \R^n\mid
(\alpha ,\lambda )=0\}\, .$$
A \textit{finite reflexion group}
is a finite subgroup $W\subset \mathrm{O}(n)$
generated by reflexions.

Let $\Delta$ be a finite set of nonzero vectors in $\mathbb{R}^n$.
Then $\Delta$ is called \textit{a root system} (in $\mathbb{R}^n$) (and
its elements are called \textit{roots}) if it satisfies
the following three conditions:

\begin{enumerate}
\item $\Delta$ generates $\mathbb{R}^n$;
\item If $\alpha \in \Delta$, then $\Delta\cap\mathbb{R}\alpha=\{\alpha,-\alpha\}$;
\item If $\alpha, \beta \in  \Delta$, then $r_{\alpha} (\beta )\in \Delta$.
\end{enumerate}

Note that sometimes it is not required that $\Delta$ generates $\mathbb{R}^n$.
The reason is that this allows one to consider subsets $\Sigma\subset \Delta$
such that $\Sigma $ is a root system in $\mathbb{R}\Sigma$ simply as root
system in $\mathbb{R}^n$, see \cite{op06} and the reference
therein for applications in analysis.

{}From now on $\Delta$ always stands for a root system in $\mathbb{R}^n$
and by $W=W(\Delta)$ we denote the reflection
group generated by the reflections $r_\alpha$, $\alpha\in\Delta$.
Then $W$ is a finite Coxeter group. Conversely, if $W$ is a finite
Coxeter group then there exists a $n\in \mathbb{N}$, and a root
system $\Delta\subset \mathbb{R}^n$,
such that $W\simeq W (\Delta)$, cf.\cite{bourLie4-6}, Chapter 1, pp 14 and 17.

Recall that a \textit{total ordering} on  a real vector space $V$ is a
transitive relation on $V$ (denoted $<$) satisfying the following axioms:

\begin{enumerate}
\item For each pair $\mu,\nu\in V$, exactly one of $\mu<\nu,\mu=\nu,\nu < \mu$ holds.
\item Let $\mu,\nu,\eta\in V$. If $\mu<\nu$, then $\mu+\eta<\nu+\eta$.
\item If $\mu<\nu$ and $c$ is a nonzero real number, then $c\mu <c\nu$ if $c>0$
and $c\nu <c\mu$ if $c<0$.
\end{enumerate}

We write $\mu >\nu$ if $\nu <\mu$.
Given such a total ordering, we say that $\nu\in V$ is \it positive \rm if $0<\nu$
and negative if $\nu <0$. Given a total ordering $<$ on $\mathbb{R}^n$
and a set of roots $\Delta\subset \mathbb{R}^n$ we
set $\Delta^{+}=\{\nu\in \Delta|0<\nu\}$. The elements in $\Delta^+$
are the \textit{positive roots}. We note that by (1) above and the fact that
$0\not\in \Delta$ it follows that $\Delta =\Delta^+\dot{\cup}(-\Delta^+)$, where
$\dot{\cup}$ stands for disjoint union.

A subset $\Pi$ of $\Delta$ is a \it simple system \rm if $\Pi$ is
a vector space basis for $\mathbb{R}^n$, and each element of $\Delta$
is a linear combination of elements of $\Pi$ with all coefficients having the same sign.
It is easy to see that if $\Pi$ is a simple system, then $w\Pi$ is also a simple system,
for any $w \in W$.
\begin{thm} Every positive system contains a unique simple system $\Pi$. Furthermore,
the Coxeter group $W(\Delta )$ is generated by the reflections $r_\alpha$, $\alpha\in \Pi$.
\end{thm}

We now  describe the construction of a fundamental domain for
the action of the Coxeter group $W$ on $\R^n$.

\begin{de} Let $G$ be a discrete group acting on $\R^n$. A closed
subset $D$ of $\R^n$ is called a fundamental domain for $G$, if
$D$ is a $G$-tile, i.e.,
$$\R^n=\bigcup_{g\in G}gD ,$$
and $gD\cap hD$ has measure zero for all $g,h \in G$, $g\not= h$.
\end{de}

\begin{de} A subset $C$ of a vector space $V$ is a cone if
$\lambda C\subseteq C$, for any real $\lambda>0$.
\end{de}
\begin{de} A subset $C$ of a vector space $V$ is
convex if for any vectors $u,v\in C$, the
vector $(1-t)u +tv$ is also in $C$ for all $t\in[0,1]$.
\end{de}
\begin{thm} Let $W=W(\Delta )$ and $\Pi\subset \Delta^+$ a simple
system of roots. Then the convex cone
$$C(\Pi)=\{\lambda\in \R^n\mid
\, (\lambda ,\alpha )\ge 0, \,\,\forall \alpha \in \Pi\}
=\{\lambda\in \R^n\mid
\, (\lambda ,\alpha )\ge 0,,\,\ \forall \alpha \in \Delta^+\}$$
is a fundamental domain for the action of $W$ on $\R^n$.
\end{thm}

Note that if we replace $\Pi$ by $w\Pi$, with $w\in W$, then
the corresponding cone is $wC$, i.e., $C(w\Pi)=w C(\Pi)$.
The open convex cones $C(\Pi)^o$ are called \textit{chambers}
and they are the connected components of the complement of
$\bigcup_{\alpha\in\Pi}H_\alpha$ in $\R^n$.
Given a chamber $C$ associated with a simple system $\Pi$, its
\textit{walls} are defined to be the
hyperplanes $H_\alpha$,  $\alpha\in\Pi$.
The angle between any two walls is an angle of the form
$\pi/k$, for some positive integer $k>1$.

We apply now the results of Section \ref{section_2} to this situation. Let
$\Delta\subset \R^n$ be a root system and $\Pi =\{\alpha_1,\ldots ,\alpha_n\}$ be
a system of simple roots. Let $\Pi^*=\{\alpha_1^*,\ldots ,\alpha_n^*\}$
be the corresponding dual basis, i.e., $(\alpha_i,\alpha_j^*)=\delta_{ij}$
for $i,j=1,\ldots ,n$. Then
$$C(\Pi )=\{t_1\alpha_1^*+\ldots +t_n\alpha_n^*\mid
t_j\ge 0,\, j=1,\ldots ,n\}\, .$$
Let $a_j>1$, $j=1,\ldots ,n$ and let $A$ be a matrix such
that $A(\alpha_j^*)=a_j\alpha_j^*$, $B=A^T$, $\cD_1=\{B^k\mid k\in \Z\}$,
and $\cD=\cD_1W$. Then $\cD_1 (C(\Pi))=C(\Pi)$.
Let $\Gamma$ be a full rank
lattice in $\R^n$ of the form $G\Z$ with $\det G=1$. We finish this
section with the following Theorem. The construction will be discussed
in some more details in Section \ref{section_5}

\begin{thm} With the notation above there exists a $(\cD,\Gamma)$ wavelet set
$\Omega\subset C(\Pi )$.
\end{thm}

\begin{proof} This follows from Theorem \ref{th-2.2}
and Theorem \ref{th-exwave}.
\end{proof}

We give now two example of Coxeter groups and we will come back to them in Section \ref{section_3}
when we will be able to describe the construction of multivavelets associated with MRA's.

\begin{Example} \label{Dihedral Group} Let $\mathbb{R}^2$ be the Euclidian plane,
and let $\cD_m$ be the dihedral group of order $2m$,
consisting of the orthogonal transformations which preserve a regular $m$-sided
polygon centered at the origin.

$\cD_m$ contains $m$ rotations through multiples of $2\pi/m$, and $m$ reflections
about the diagonals of the polygon. By 'diagonal', we mean
a line joining two vertices or the midpoints of opposite sides if $m$ is even, or
joining a vertex to the midpoint of
the opposite side if $m$ is odd.

The group $\cD_m$ is actually generated by reflections, since a rotation through
$2\pi/m$ is a product of two
reflections relative to a pair of adjacent diagonals which meet at an angle of
$\theta=\pi/m$, see Figure \ref{Dm}.
\begin{figure}[h]
\begin{center}
\epsfig{scale=0.5,angle=0,file=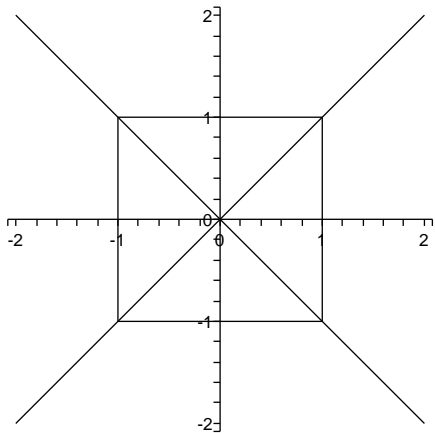}
\caption{  The dihedral group $\mathcal{D}_4$ }
\label{Dm}
\end{center}
\end{figure}
\end{Example}

The three dimensional case is more interesting. Let $a,b,c$ be three linearly
independent vectors
such that the corresponding reflections lie in a finite group. That is only
possible if $\sphericalangle(a,b)$,
$\sphericalangle(a,c)$, $\sphericalangle(b,c)$ are rational multiple of
$\pi$. This can be obtain by choosing $\sphericalangle(a,b)$
to be an arbitrary multiple of $\pi$ and then choosing $c$ such that
$\sphericalangle(a,c)=\sphericalangle(b,c)=\pi/2$.
In that case, the group generated by $r_a$ and $r_b$, $<r_a,r_b>$, is a dihedral group and
$<r_a,r_b,r_c>$ is the direct product of the dihedral group $<r_a,r_b>$ and the cyclic group of
order 2 generated by $r_c$.

Except these direct products, there are only three $3$-dimensional Euclidian reflection groups,
the groups of symmetries of a regular tetrahedron, a cube, and a regular dodecahedron.
\begin{Example} \label{Tetrahedron}
 For each tetrahedron centered, there is a 'dual' tetrahedron which
is congruent to the given one and has the property that each edge of the given tetrahedron is
perpendicularly bisected by an edge of the dual. Together, the vertices
of the two tetrahedra give the vertices of
a cube. Let $a$ and $c$ be the position vectors of the midpoints of a
pair of parallel but not opposite edges $e_1$
 and $e_2$. Let $b$ be the position vector of the midpoint of one of the edges on the
opposite face to that determined by $e_1$ and $e_2$, which are not parallel to $e_1$ and $e_2$.
Then we have the following:
\begin{displaymath}
\begin{array}{ll}
\sphericalangle(a,b)=2\pi/3, &  r_ar_b\, \, \hbox{has order 3};\\
\sphericalangle(a,c)=\pi/2, &  r_ar_c\, \, \hbox{has order 2};\\
\sphericalangle(b,c)=2\pi/3, &  r_br_c\,\,  \hbox{has order 3};\\
\end{array}
\end{displaymath}

$r_a,r_b,r_c$ are symmetries of the tetrahedron,
and so these three reflections generate the group of all symmetries of the
tetrahedron, which is just Sym(4).
\end{Example}

\section{Multiwavelets Associated with Multiresolution Analysis}\label{section_3}
\noindent
This section and the next one contain our main results. We start by constructing
scaling sets and by taking translates we obtain
MRA multi-wavelets which are completely supported in a fundamental domain for the
action of a Coxeter group on $\mathbb{R}^n$.
Let $A$ be an expansive matrix and $\mathcal{T}$ be a full rank lattice.
\begin {de} A multiresolution analysis on $\mathbb{R}^n$ is a sequence
of subspaces $\{ V_j\}_{j\in\mathbb{Z}}$ of functions in
$L^2(\mathbb{R}^n)$ satisfying the following properties:
\begin{deflist}
\item[i)]For all $j\in\mathbb{Z}$, $V_j\subseteq V_{j+1}$
\item[ii)]If $f(\cdot)\in V_j$, then $f(A\cdot)\in V_{j+1}$
\item[iii)] $\bigcap_{j\in\mathbb{Z}}{V_j}=\{0\}$
\item[iv)] $\overline{\bigcup_{j\in\mathbb{Z}}{V_j}}=L^2(\mathbb{R}^n)$
\item[v)] There exists a function $\phi\in L^2(\mathbb{R}^n)$ such that
$\{\phi(\cdot+t)|t\in\cT\}$ is an orthonormal basis for $V_0$.
\end{deflist}
\end {de}
The function $\phi$ is called a \textit{ scaling function}.
One can allow more than one scaling function, say $m$, \rm and then the
MRA has multiplicity \it m. \rm

If we change condition \it iv) \rm in the definition of the multiresolution analysis into

\begin{deflist}
\item[iv')] $\overline{\bigcup_{j\in\mathcal{T}}{V_j}}=L^2_{M}(\R^n) $, for some subset
$M\subseteq\R^n$,
\end{deflist}
then we get a
\textit{subspace multiresolution analysis}, SMRA.

We explain now how to construct a wavelet or
multi-wavelets from a MRA.
Let $W_0$ be the orthogonal complement of $V_0$ in $V_1$, that is, $V_1=V_0 \oplus W_0$.
In general, let $W_i=V_{i+1} \ominus V_{i}$, for each $j\in \mathbb{Z}$, $$V_j=
\bigoplus\limits_{l=-\infty}^{j}W_l$$ and so $$L^2(\mathbb{R}^n)=
\bigoplus\limits_{l=-\infty}^{\infty}W_l.$$
If there exists a function $\psi\in W_0$ such that $\{\psi(\cdot+t)| t\in\mathcal{T} \}$ is an
orthonormal basis for $W_0$, then $\{\psi_{j,t}| t\in\mathcal{T} \}$ is an
orthonormal basis for $W_j$, and $\{\psi_{j,t}| t\in\mathcal{T}, j\in\mathbb{Z}\}$ is an
orthonormal basis for $L^2(\mathbb{R}^{n})$, which means that
$\psi$ is an orthonormal wavelet associated with the given MRA..

Recall that $\mathcal{T}$ is a full rank lattice and $A$ an expansive matrix. Set $B=A^T$.
Suppose $B\mathcal{T}\subseteq \mathcal{T}$. Let $\mathcal{T}/{B\mathcal{T}}$ be
the quotient group, where we identify its elements with their
representative vectors in $\mathbb{R}^n$, $v_0,v_1,...v_{q-1}$, where $q=|\det(B)|$.
\begin{lem}
Let $K$ be a $\mathcal{T}$-tile such that $B^{-1}K \subset K$.
Let $$K_i=(B^{-1}K+B^{-1}v_i+\mathcal{T})\bigcap K.$$
Then
\begin{deflist}
\item[i)]$K=\bigcup_{i=0}^{q-1}K_i$ up to measure zero and $K_i\bigcap K_j=0$ up
to measure zero for $i\ne j$.
\item[ii)]$BK_i\sim_{\mathcal{T}} K$.
\end{deflist}
\begin{proof}Let $x\in {K_i\cap K_j}$.
Then there exist $u_1,u_2\in K$, $v_i,v_j\in\mathcal{T}/{B\mathcal{T}}$ and
$t_1,t_2 \in \mathcal{T}$, such that
$$x = B^{-1}u_1+B^{-1}v_i+t_1=B^{-1}u_2+B^{-1}v_j+t_2,$$ and so $u_1-u_2\in \mathcal{T}$.
But $K$ is a $\mathcal{T}$ tile, so $u_1=u_2$, and then $v_i-v_j=B(t_1-t_2)$. Thus $v_i=v_j$
and so $i=j.$
Let the notations be as before and let $$K_{i,t}=B^{-1}K\bigcap (K-B^{-1}v_i-t).$$
Then $\bigcup_{t\in\mathcal{T}}(K_{i,t}+B^{-1}v_i+t)=K_i$.

Since $K$ is a $\mathcal{T}$ tile, it follows that $K-B^{-1}v_i$ is also a
$\mathcal{T}$ tile. Thus $K_{i,t}$ are measurewise disjoint and
$\bigcup_{t\in\mathcal{T}}K_{i,t}=B^{-1}K$.

By definition, $K_i\subseteq K$, and $$|K_i|=\Sigma_{t\in\mathcal{T}}|K_{i,t}|=
|B^{-1}K|=q^{-1}|K|$$ for $i=0,...,q-1$, and $|K_i\bigcap K_j|=0$ for $i\ne j$.
Therefore  $$K=\bigcup_{i=0}^{q-1}K_i$$ up to measure zero.
Moreover, $$B(K_i)=\bigcup_{t\in\mathcal{T}}B(K_{i,t}+B^{-1}v_i+t)=K+v_i+Bt$$ and thus
$$B(K_i)\sim_{\cT}K.$$
\end{proof}
\end{lem}
Let $K\subset\R^n$ be a measurable set. Set $V_0=L^2_K(\R^n)$ and
$V_j=\{f(A^j\cdot)|f(\cdot)\in V_0\}$.

\begin{de}A set $K\subset\R^n$, $|K|=1$, is a scaling set, if the sequence $\{V_j\}$ described
above
is a multiresolution analysis with scaling function $\phi=\cF^{-1}\chi_{K}$.

\end{de}
\begin{thm}
A subset $K\subset\mathbb{R}^n$ is a scaling set if and only if
$B^{-1}K\subseteq K$ and $K$ is a $\mathcal{T}-tile$.
\end{thm}
\begin{proof} Suppose first that $K\subset\mathbb{R}^n$ is a scaling set.
Then $\phi=\cF^{-1}{\chi_K}$ is a scaling function and so
$\{\phi_{0,t}\}_{t\in\mathcal{T}}$
is an orthonormal basis for $V_0$. This implies that
$\{\hat{\phi}_{0,t}\}_{t\in\mathcal{T}}$
is an orthonormal basis for $\hat{V_0}$.
For $t\in\R^n$ let $e_t(x)=e^{2\pi (x,t)}$.
It follows that $\hat{V_0}$ has a orthonormal basis
of the form $\{e_t \chi_K\}_{t\in\mathcal{T}}$.
{}From this we get two things. The first is that $(K,\mathcal{T})$ is a
spectral pair and thus $K$ is a $\mathcal{T}$-tile. The second is that
$\hat{V_0}=L^2(K)$ and by the SMRA structure, we get that
$\hat{V}_{-1}\subset\hat{V}_{0}$
which implies that
$B^{-1}K\subset K$.

Assume now that $B^{-1}K\subset K$ and that $K$ is a $\mathcal{T}-tile$.

Set $\hat{V_0}=L^2(K)$ and $\hat{V_j}=L^2(B^jK)$.
Since $B^{-1}K\subset K $, it follows that $\hat{V}_{j}\subset\hat{V}_{j+1}$.
The other conditions are easy to verify. Thus
$K\subset\mathbb{R}^n$ is a scaling set.
\end{proof}

The next theorem gives in a constructive way, the existence of SMRA wavelets.
\begin{thm}
If $K\subset\mathbb{R}^n$ is a scaling set, then
\begin{deflist}
\item[i)]$\hat{V_0}=L^2(K)$ and $\hat{V_j}=L^2(B^jK)$
\item[ii)]$\{\psi^{i}=\check{\chi}_{\Omega_i}\}_{i=1}^{q-1}$ is a
SMRA multiwavelet, where $\Omega_i=BK_i$.
\end{deflist}
\end{thm}
\begin{proof}i) It follows from the theorem above.

ii) By lemma 1, $\Omega_i\sim K$. This implies that $(\Omega_i,\mathcal{T})$
is a spectral pair. Then $$\{\hat{\psi}_{o,t}^i\}_{t\in\mathcal{T}}$$ is an
orthonormal basis for $L^2(\Omega_i)$.

Set $\widehat{W}_{0,i}=L^2(\Omega_i)$ for $i=1,...,q-1$.
By construction, $$BK= K \cup {\bigcup_{i=1}^{q-1}\Omega_i}.$$
Therefore $$\widehat{V}_1=L^2(BK)=L^2(K)\oplus \bigoplus_{i=1}^{q-1}L^2(\Omega_i)
=\widehat{V}_0\oplus\widehat{W}_{0,1}\oplus...\oplus\widehat{W}_{0,q-1}.$$
So
$$V_1=V_0\oplus\bigoplus_{i=1}^{q-1}W_{0,i}$$
and for any $j\in \mathcal{Z}$,
we have
$$V_{j+1}=V_j\oplus\bigoplus_{i=1}^{q-1}W_{j,i}\, .$$
Thus,
$\{\psi^{i}=\cF^{-1}{\chi}_{\Omega_i}\}_{i=1}^{q-1}$ is a SMRA multi-wavelet.
\end{proof}

\section{Multiresolution and Coxeter Groups}\label{section_4}
\noindent
We now go back to the notation in Section \ref{section_1}.
In particular $\Delta$ is a root system
in $\R^n$, $\Pi$ a system of simple roots,
$W=<r_{\alpha_i}|\alpha_i\in\Pi>$ the corresponding
a finite Coxeter group and $C=C (\Pi)$ be the corresponding
positive cone which is a  fundamental domain for the
action of $W$ on $\R^n$. Let $A$ be such that with $B=A^T$
we have
$B(\alpha_j^*)=a_j \alpha_j^*$, $a_j>1$, $j=1,\ldots ,n$.
The more general case $B(\alpha_j)=a_j\alpha_{\pi (j)}^*$ where,
$\pi :\{1,\ldots ,n\}\to \{1,\ldots ,n\}$ is
a permutation, is handled in similar way. Let
$$P=\{\sum_{i=1}^n t_i\alpha_i^*|0<t_i\leq s_i\},$$ where $s_i$ are
such that $|P|=1$. Note that $P$ is a $n$ dimensional
parallelepiped and a $B \mathbb{Z}^n $ tile.

Indeed, if $z\in \Z^n$, then

\begin{eqnarray*}
P+Rz & = & \{\sum_{i=1}^n t_i\alpha_i^*+\sum_{i=1}^nn_i\alpha_i^*|0<t_i\leq s_i,\,n_i\in\Z\}\\
& = & \{\sum_{i=1}^n (t_i+n_i)\alpha_i^*|0<t_i\leq s_i,\,n_i\in\Z\}
\end{eqnarray*}

so
$$|P\cap(P+Rz)|=0.$$
We have
$$B^{-1}P=\{\sum_{i=1}^n a_i^{-1}t_i\alpha_i^*|o<t_i\leq s_i\}\subset P\, ,$$
since $0<a_i^{-1}t_i\leq s_i$.
Moreover,
$$(\sum_{i=1}^n t_i\alpha_i^*,\alpha_m)=\sum_{i=1}^n t_i(\alpha_i^*,\alpha_m )=t_m>0\, .$$
Thus, $B^{-1}P\subset P\subset C (\Pi) $.

\begin{thm}Let $P,B$ be as above, $P_i=B^{-1}P+B^{-1}v_i$ and  $\Omega_i=BP_i$. Then
$\{\psi^{i}=\check{\chi}_{\Omega_i}\}_{i=1}^{q-1}$ is a
SMRA multiwavelet associated to the multiresolution $V_j=L^2_{B^jP}(\R^n)$.
\end{thm}
\begin{proof}
As shown above, $P$ is a $R\mathbb{Z}^n $-tile and $B^{-1}P \subset P$ and so $P$ is a scaling set.
Thus, by theorem 4.8, $\{\psi^{i}=\check{\chi}_{\Omega_i}\}_{i=1}^{q-1}$ is a SMRA multi-wavelet.
\end{proof}

\begin{figure}
\begin{center}
\epsfig{scale=0.8,angle=0,file=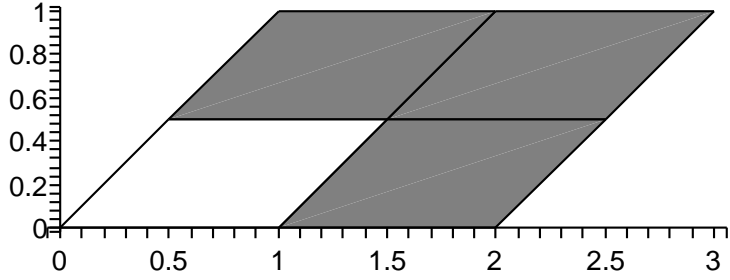}
\caption{SMRA wavelet sets in $\R^2$}
\label{rot2}
\end{center}
\end{figure}
\begin{Example}
Let the group $\mathcal{D} =  \{R^{k}_{\frac{2\pi}{m}}\}_{k=0}^{m-1}$ act on $\R^2$ and let
$D=\{ t_1(1,0)+t_2(\cot2\pi/m,1), \,0< t_{1,2}\}$ be the fundamental domain of this action.
Let $B=2\id_2$ and $P=\{ t_1(1,0)+t_2(\cot2\pi/m,1), \, t_{1,2}\in[0,1]\}$.
Let $$\Omega_1=P+(1,0),$$ $$\Omega_2=P+(1+\cot2\pi/m,1),$$$$\Omega_3=P+(\cot2\pi/m,1).$$
Then $\{\psi^{i}=\check{\chi}_{\Omega_i}\}_{i=1}^{3}$ is a SMRA multiwavelet, see Figure \ref{rot2}.
\end{Example}

\begin{Example}
Let $\Omega=<r_a,r_b,r_c>$ be a Coxeter group, where $a,b,c$ are as described in  Example \ref{Tetrahedron}.
Then the fundamental domain for the action of $\Omega$ on $\R^3$ is
$$D=\{ t_aa^*+t_bb^*+t_cc^*|0<t_a,t_b,t_c\}.$$
Let $P=\{ t_aa^*+t_bb^*+t_cc^*|0<t_a<s_a,\,0<t_b<s_b,\,0<t_c<s_c\}$, such that $|P|=1$.

Let $B=2\id_3$. Then $\det(B)=2^3=8$ and so there are 7 MRA wavelet sets, see Figure \ref{threedim}.

\begin{figure}[h!]
\begin{center}
\epsfig{scale=0.9,angle=0,file=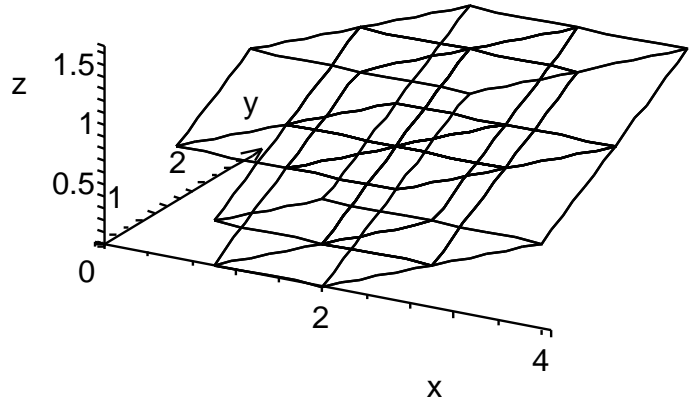}
\caption{SMRA wavelet sets in $\R^3$}
\label{threedim}
\end{center}
\end{figure}
\end{Example}

\begin{rem} We would like to remark that the
constructions in this section  can easily be
reduced to the case $\alpha_i=e_i$, where $e_i$ is the
standard basis for $\R^n$. This is done by
using the linear map $\alpha_j\mapsto e_j$.
 \end{rem}

\section{Wavelet Sets and Coxeter Groups}\label{section_5}
\noindent
We now explain how the construction of wavelet sets mentioned in
Section \ref{section_2} can be done.
\begin{thm}Let $P$ be as above and let $F=BP\setminus P$.
Define
$$W_{1,1}=(P\setminus B^{-1}P)+\alpha_i^{*}$$
$$W_{2,1}= B^{-2}[F\setminus (P+\alpha_i^{*})]$$
$$\ldots $$
$$W_{1,n}= [(B^{-n+1}P\setminus B^{-n}P)\setminus W_{2,n-1}]+\alpha_i^{*}$$
$$W_{2,n}= B^{-n-1}\{[(B^{-n+1}P\setminus B^{-n}P)+\alpha_i^{*}]\setminus W_{1,n}\}.$$
Then $$P=\bigcup_n^\infty W_{2,n}\bigcup\bigcup_n^\infty (W_{1,n}-\alpha_i^{*})$$
$$F=\bigcup_n^\infty W_{1,n}\bigcup\bigcup_n^\infty B^{n+1}W_{1,n}.$$
Moreover, if we let$$W=\bigcup_{j=1,2}\bigcup_n^\infty W_{j,n},$$
then $W$ is a wavelet set.

\end{thm}
\begin{proof}
We have shown above that $P$ is a $R\mathbb{Z}^n $-tile. On the other hand,
$$BF\cap F=\emptyset,$$ and
$$ \bigcup_{n\in\Z}B^nF=D,$$
so $F$ is a multiplicative tiling.
By definition,
$$W\sim_{R\mathbb{Z}^n}P,$$
and $$W\sim_{B}F.$$
Thus, $W$ is a wavelet set.
\end{proof}

\section{Coxeter Groups and Sampling}
\noindent
In this final section we discuss how the results in Section \ref{section_4} are
related to sampling theory. First we note that
$$L^2(\R^n)=\bigoplus_{w\in W}L^2_{C(w\Pi )}(\R^n)\, .$$
Thus each $f\in L^2(\R^n)$ can be decomposed as
\begin{equation}\label{eq-6.1}
f=\sum_{w\in W}f_w,
\end{equation}
where $f_w=\cF^{-1}(\chi_{C(w\Pi)}\cF(f))$ is the orthogonal
projection of $f$ onto $L^2_{C(w\Pi )}(\R^n)$. The function $f_w$ contains the
exact frequency information of $f$ in the direction of the cone
$C(w\Pi )$.

We note that $C(w\Pi )^o$ is an open convex cone such that $C(w\Pi)^o\cap -C(w\Pi )^o=
\emptyset$, thus the dual cone $C_w\not=\emptyset$ is an open convex cone. In fact
\begin{eqnarray*}
C_w:=[C(w\Pi )^o]^*&=&\{y\in \R^n\mid (\forall \lambda \in
C(w\Pi)^o)\, (y,\lambda )>0\}\\
&=&w\{\sum_{i=1}^nt_j\alpha_j\mid t_j>0,\, j=1,\ldots ,n\}\\
&=&w[C(\Pi )^o]^*\not=\{0\}
\end{eqnarray*}
is an open convex cone, and the function
$f_w$ extends to a holomorphic function on the tube domain
$$T(C_w)= \R^n+iC_w\, .$$
This holomorphic extension is given by
\begin{eqnarray*}
F_w(x+iy)&=&\int_{C (w\Pi )}\cF (f)(\lambda )e^{2\pi i (x+iy,\lambda )}\, d\lambda\\
&=&\int_{C(w\Pi )}\cF (f)(\lambda )e^{-2\pi (y,\lambda )}e^{2\pi i (x,\lambda )}\, d\lambda
\end{eqnarray*}
and
\begin{equation}\label{eq-6.2}
f_w(x)=\lim_{y\to 0}F_w(x+iy )
\end{equation}
where the limit is taken in $L^2(\R^n)$, see \cite{SW71} for details.
Thus $f_w$ is in the Hardy space $H^2(T(C_w))$ and the equations
(\ref{eq-6.1}) and (\ref{eq-6.2}) give us a
decomposition of $f$ as $L^2$-limit of $\# W$-holomorphic functions.

We can then sample an approximate version of each $f_w$ (and then $f$) using the
following simple version of the Whittaker-Shannon-Kotel'nikov sampling
theorem \cite{Sh48,Za93}. We include a short proof using
the idea of spectral sets which we have not seen elsewhere in the
literature, even if it is well known. Note that if the Fourier transform of $f$ is supported
in a set of finite measure, then $f$ is continuous so $f(x)$ is well
defined for all $x\in\R^n$.
\begin{thm}(WSK-sampling theorem for spectral sets)\label{wsk}
Let $P\subset \R^n$ be measurable, $0<|P|<\infty$ and such that
there exists a discrete set $\Gamma\subset \R^n$, such that the functions
$\{e_\gamma\chi_P\}_{\gamma\in\Gamma}$ form a orthogonal basis for $L^2(P)$, i.e.,
$P$ is similar to a spectral set. Let
$$\varphi =|P|^{-1/2}\cF^{-1}\chi_P\, .$$
Then we have that for all $f\in L^2_P (\R^n)$,
$$(\cF(f),e_\gamma )=f(-\gamma)$$
and
$$f(x)=|P|^{-1}\sum_{\gamma\in\Gamma}f(-\gamma )\varphi (x+\gamma)=
|P|^{-1}\sum_{\gamma\in -\Gamma}f(\gamma )\varphi (x-\gamma)\, .$$
Furthermore, if $B\in\GL (n,\R)$ and $L=(B^{-1})^T$. Then
\begin{eqnarray*}
f(x)&=&\frac{1}{|P|}\sum_{\gamma\in\Gamma}f(-L\gamma )\varphi (B^Tx+\gamma)\\
&=&\frac{1}{|P|}\sum_{\gamma\in\Gamma}f(-L\gamma )\varphi (B^T(x+L\gamma))\\
&=&|P|^{-1}\sum_{\gamma\in L\Gamma}f(-\gamma )\varphi\circ B^T (x+\gamma)
\end{eqnarray*}
for all $f\in L^2_{BP}(\R^n)$.
\end{thm}
\begin{proof} We have
$$(\cF(f),e_\gamma\chi_P)= \int_{\R^n}\cF (f)(\lambda )e^{-2\pi (\lambda ,\gamma)}\, d\lambda
=f(-\gamma )\, .$$
Furthermore, as $\{|P|^{-1/2}\, e_\gamma \chi_P\}_{\gamma\in\Gamma}$ is an orthonormal
basis for $L^2(P )$,
\begin{eqnarray*}
\cF(f)(\lambda )&=&|P|^{-1}\sum_{\gamma\in \Gamma} (\cF (f),e_\gamma \chi_P)
(e_\gamma\chi_P) (\lambda )\\
&=&|P|^{-1}\sum_{\gamma\in \Gamma}
 f(-\gamma )(e_\gamma\chi_P) (\lambda )\, .
\end{eqnarray*}
Hence
\begin{eqnarray*}
f(x)&=&\sum_{\gamma\in \Gamma}f(-\gamma )\cF(e_\gamma\chi_P) (x )\\
&=&\sum_{\gamma\in \Gamma}f(-\gamma )\varphi (x+\gamma )\, .
\end{eqnarray*}
Let now $f\in L^2_{BP}(\R^n)$.
As $\cF(f\circ L)(\lambda )=|\det L|^{-1}\cF (f)(B\lambda )$
it follows that $f\circ L\in L^2_P(\R^n)$. Hence
\begin{eqnarray*}
f(x)&=&(f\circ L)(B^Tx)\\
&=&|P|^{-1}\sum_{\gamma\in\Gamma}(f\circ L)(-\gamma )\varphi (B^Tx+\gamma)\\
&=&|P|^{-1}\sum_{\gamma\in\Gamma}f(-L \gamma )|\det B|\varphi (B^T(x+L\gamma))
\end{eqnarray*}
and the claim follows.
\end{proof}

We now go back to the situation in Section \ref{section_4}.
 Let $\Pi =\{\alpha_1,\ldots ,\alpha_n\}$ be a system
of simple roots and $\Pi^*=\{\alpha_1^*,\ldots ,\alpha_n^*\}$ its dual basis.
If we let $P=\{\sum_{j=1}^nt_j \alpha_j^*\mid 0<t_j\le s_j\}$, where the numbers $s_j$ are
chosen so that $|P|=1$, then $P$ is a spectral set and
$\{e_\gamma\chi_P\}_\gamma$ is a orthonormal basis of $L^2(P )$ where
$\Gamma = \Z \Pi$. Let $B(\alpha_j^*)=a_j \alpha_j^*$ be as before
and $V_j=L^2_{B^jP}(\R^n)$. By Theorem \ref{wsk} we therefore get:
\begin{thm} Let $\varphi =\cF^{-1}\chi_P$ and $L=(B^{-1})^T$.
Then, if $f\in V_j=L^2_{B^jP}(\R^n)$
we have
\begin{eqnarray*}
f(x)&=&\sum_{\gamma\in\Gamma}f(-L^j\gamma )\varphi ((B^j)^Tx+\gamma)\\
&=&\sum_{\gamma\in\Gamma}f(-L^j\gamma )\varphi ((B^j)^T(x+L^j\gamma))\\
&=& \sum_{\gamma\in L^j\Gamma}f(-\gamma )\varphi\circ B^T (x+\gamma)\, .
\end{eqnarray*}
\end{thm}
\bibliographystyle{plain}

\def\cprime{$'$} \def\polhk#1{\setbox0=\hbox{#1}{\ooalign{\hidewidth
  \lower1.5ex\hbox{`}\hidewidth\crcr\unhbox0}}}

\end{document}